\documentclass{amsart}
\usepackage{amsfonts}
\usepackage{amssymb}
\usepackage{amsmath}
\usepackage{amsrefs}
\usepackage{mathrsfs}

\usepackage{color}
\usepackage[dvipsnames]{xcolor}
\usepackage{tikz}
\usepackage{setspace}
\usepackage{multicol}
\usetikzlibrary{calc,intersections,decorations.pathreplacing}
\usepackage{ragged2e}
\usepackage[linktocpage=true]{hyperref}
\usepackage{pgfplots}

\usetikzlibrary{arrows.meta}
\usetikzlibrary{arrows}

\newtheorem{thm}{Theorem}[section]
\newtheorem{cor}[thm]{Corollary}
\newtheorem{lem}[thm]{Lemma}
\newtheorem{prop}[thm]{Proposition}
\newtheorem{prob}[thm]{Problem}

\theoremstyle{definition}

\newtheorem{prop-def}[thm]{Proposition--Definition}
\newtheorem{claim}[thm]{Claim}

\theoremstyle{remark}
\newtheorem{rmk}[thm]{Remark}

\usetikzlibrary{fit}

\hypersetup{ colorlinks=true, citecolor=cyan, linkcolor=blue, urlcolor=blue}

\begin{document}

\title[Rational curves on del Pezzo manifolds]{Rational curves on del Pezzo manifolds}

\author[A. Zahariuc]{Adrian Zahariuc}
\address{Department of Mathematics, Harvard University, Cambridge, MA 02138, USA}
\curraddr{Department of Mathematics, UC Davis, One Shields Ave, Davis, CA 95616}
\email{azahariuc@math.ucdavis.edu}
\thanks{Partially supported by National Science Foundation grant DMS-1308244, "Nonlinear Analysis on Sympletic, Complex Manifolds, General Relativity, and Graphs".}





\begin{abstract}
We exploit an elementary specialization technique to study rational curves on Fano varieties of index one less than their dimension, known as del Pezzo manifolds. First, we study the splitting type of the normal bundles of the rational curves. Second, we prove a simple formula relating the number of rational curves passing through a suitable number of points in the case of threefolds and the analogous invariants for del Pezzo surfaces.
\end{abstract}

\maketitle

\vspace*{6pt}\tableofcontents

\section{Introduction}

\subsection{Statement of the main results} 
The index of a projective Fano variety $X$ is the largest integer by which $-K_X$ is divisible in $\mathrm{Pic}(X)$. Let $n = \dim_{\mathbf C} X$. It can be proved that the index is at most $n+1$, and moreover, it is equal to $n+1$ only for projective spaces, respectively to $n$ only for smooth quadric hypersurfaces. 

We will be concerned with the study of the rational curves on $X$. Hypersurfaces of degree one or two are homogeneous convex varieties, which are quite easily understood from the point of view of quantum cohomology, so one might guess that Fano varieties of index $n-1$ are the simplest varieties for which the behavior of rational curves is harder to control. In this paper, we will use an elementary specialization technique to study some basic enumerative questions concerning rational curves on such algebraic varieties. I know of no analogue of the approach for Fano varieties of smaller or larger index.

A smooth Fano projective variety $X$ such that $-K_X$ is divisible by $n - 1$ in $\text{Pic}(X)$ is sometimes called a del Pezzo variety or del Pezzo manifold. These varieties can be thought of as higher-dimensional analogues of del Pezzo surfaces. The hyperplane class ${\mathscr O}(1)$ is uniquely defined by the property $-K_X={\mathscr O}(n-1)$. The degree of $X$ is of course $\smash{ d = \int_{X} c_1({\mathscr O}_X(1))^n}$. Del Pezzo varieties have been completely classified by Fujita and Iskovskikh ~\cite{Fu80, Fu90, Is78, Is80}. In this paper, we will mainly be concerned with the classes of Picard rank one and base point free polarization, which amounts to $d \in \{2,3,4,5,8\}$, to avoid excessive bookkeeping when $d \in \{6,7\}$, respectively some technical issues when $d=1$. 

$\bullet$ if $d=2$, then $X$ is a double cover of ${\mathbb P}^n$ ramified along a quartic;

$\bullet$ if $d=3$, then $X$ is a cubic hypersurface in ${\mathbb P}^{n+1}$;

$\bullet$ if $d=4$, then $X$ is a complete intersection of two quadrics in ${\mathbb P}^{n+2}$;

$\bullet$ if $d=5$, then $X$ is a linear section of the Grassmannian ${\mathbb G}(1,4) \subset {\mathbb P}^9$;

$\bullet$ if $d=8$, then $X$ is ${\mathbb P}^3$.
\newline In the last case, the polarizing line bundle ${\mathscr O}(1)$ is what we would normally call ${\mathscr O}_{{\mathbb P}^3}(2)$ and the index is $4$ rather than $2$.

When useful, we denote a del Pezzo manifold of degree $d \in \{2,3,4,5,8\}$ and dimension $n \geq 2$ by $\mathrm{dP}_n^d$. Furthermore, if $\mathrm{dP}_m^d$ is a smooth plane section of $\mathrm{dP}_n^d$, which is still del Pezzo by adjunction, we denote the inclusion by $j_{m,n}^d : \mathrm{dP}_m^d \to \mathrm{dP}_n^d$.

\begin{thm}\label{1.1}
Assume that $j_{2,n}^d: \mathrm{dP}_2^d \to \mathrm{dP}_n^d$ is a sufficiently general plane section of a del Pezzo manifold of degree $d \in \{2,3,4,5,8\}$. Let $\smash{ f: {\mathbb P}^1 \to \mathrm{dP}_2^d }$ be an unramified degree $e$ morphism such that $\smash{ f_*[{\mathbb P}^1] \notin c_1(\mathrm{dP}_2^d) {\mathbf Q} }$. Then
$$ {\mathscr N} \cong {\mathscr O}(e-1)^{\oplus 2} \oplus {\mathscr O}(e)^{\oplus (n-3)}, \eqno(1.1) $$
where ${\mathscr N} = {\mathscr N}_{j_{2,n}^d \circ f}$ is the normal sheaf of of $f$ relative to $\mathrm{dP}_n^d$.
\end{thm}

\noindent This statement generalizes an old result of Harris, Hulek and Eisenbud and Van de Ven ~\cite{EV81, Hu81} which says that a type $(1,k)$ rational curve on a smooth quadric surface has balanced normal bundle relative to the ambient ${\mathbb P}^3$ if $k > 1$. 

An immediate corollary is that, with few obvious exceptions which are left to the reader, for a homology class $\beta \in \mathrm{H}_2(\mathrm{dP}_n^d,\mathbf Z)$, there exists an irreducible component of $\overline{\mathcal M}_{0,0}(\mathrm{dP}_n^d,\beta)$ such that a general stable map in that component has balanced normal bundle. Irreducibility is actually known to hold for the cubic hypersurface case for $n \geq 4$ and to almost hold for $n=3$ ~\cite{CS09}. The reader may consult ~\cite{EV81, Ra07, Sh12} for other results concerning the property of the normal bundles of being generically balanced. 

\begin{thm}\label{1.2}
Let $\mathrm{dP}_3^d$ be a Fano threefold of index $2$ and degree $d \in \{2,3,4,5\}$, a homology class $\beta \in \mathrm{H}_2(\mathrm{dP}_3^d,\mathbf Z)$ and $e = \deg \beta$. Then we have
$$ \left\langle [\mathrm{pt}]^e \right\rangle^{\mathrm{dP}_3^d}_{0,\beta} =  \frac{1}{d(9-d)} \sum_{(j_{2,3}^d)_*\gamma =\beta} \left((K \cdot \gamma)^2 - \left(K^2\right)\left(\gamma^2\right) \right) \left\langle[\mathrm{pt}]^{e-1} \right\rangle^{\mathrm{dP}_2^d}_{0,\gamma} \eqno(1.2) $$
where $K$ is the canonical divisor of the del Pezzo surface $\mathrm{dP}_2^d$.
\end{thm}

\noindent The invariants above can be computed using more general enumerative methods -- in particular, the genus zero GW invariants of rational surfaces can be computed from the associativity of the quantum product ~\cite{GP98}. However, the identity itself is new and surprisingly simple. The most mysterious feature of the identity is the constant $d(9-d)$, whose existence is natural from the point of view of the proof of Theorem \ref{1.2}, but whose nice exact value seems to come down to some numerical coincidences. It would be interesting to find a clear explanation, if one exists. 

There is an analogous identity in the case of ${\mathbb P}^3$, which was found -- or rather conjectured -- by Coray ~\cite{Co83} in the early 80's (!) using a different argument. Very recently, the identity for ${\mathbb P}^3$ was proved by Brugall\'{e} and Georgieva ~\cite{BG16}, following an approach suggested by Koll\'{a}r ~\cite{Ko15}, which is similar to the one used here. Finally, we note that our numbers agree with those computed in ~\cite{Be95}, for $d \in \{3,4\}$ and $e \leq 3$.

\subsection{The approach by specialization} For the reader's convenience, we outline the idea of the proof in the case of a cubic hypersurface $X$ in ${\mathbb P}^{n+1}$. The generalization to other classes is completely straightforward. Ignoring any technical or compactification-related subtleties for now, let ${\mathcal C} \to {\mathcal M}$ be the moduli space of degree $e$ rational curves with its universal family. Consider the incidence correspondence $I$ between degree $e$ rational curves and $e$-tuples of points on $X$ such that all $e$ points lie on the curve, as shown in the diagram below.
\begin{center}
\begin{tikzpicture}
\matrix [column sep  = 5mm, row sep = 8mm] {
	& & & & \node (c) {$I$}; & \\
	\node (ww) {$E^e$}; & & & \node (sw) {$X^e$}; &
	& \node (se) {${\mathcal M}$}; \\
};
\draw[right hook-{Stealth[length=2mm]}, thin] (ww) -- (sw);
\draw[-{Stealth[length=2mm]}, thin] (c) -- (sw);
\draw[-{Stealth[length=2mm]}, thin] (c) -- (se);

\node at (3.3,0.35) {$=\underbrace{{\mathcal C} \times_{\mathcal M} {\mathcal C} \times_{\mathcal M} ... \times_{\mathcal M} {\mathcal C}}_{e \text{ times}}$};
\end{tikzpicture}
\end{center}  
The idea is to specialize the $e$ points to $e$ distinct points on a smooth plane cubic curve $E$ on $X$ obtained by cutting $X$ with a $2$-plane. We can understand the fiber of $I \to X^e$ over such a special $e$-tuple thanks to the following observations: 

(1) essentially by B\'{e}zout's theorem, the curves containing the $e$ special points are forced to lie on a cubic surface obtained by cutting $X$ with a $3$-plane; 

(2) the divisor class of the curve on the cubic surface restricts on $E$ to the divisor class of the $e$ points; and 

(3) the fact that $E$ has positive genus ensures that components of degenerate rational curves can't fall inside $E$, so to speak, unless they are contracted. The latter possibility turns out not to interfere with our affairs. 

The failure of property (3) for curves of positive genera is the main impediment to extending this method to arbitrary genus, but it is conceivable that excess intersection could make some headway. Finally,  if possible, it would be very interesting to take the "direct product" of ~\cite{BG16} and this paper, in order to compute the Gromov-Witten-Welschinger invariants of Fano threefolds of index $2$.  

\subsection{The generating functions} We conclude the introduction with the simple observation that the formula in Theorem \ref{1.2} can be naturally expressed using generating functions. Consider the rings of formal power series on $\mathrm{H}^2(\mathrm{dP}_n^d,{\mathbf R})$ with an additional formal variable whose sole purpose is to ensure convergence
$$ F_n^d[[t]] = {\mathbf R}[[t]] \otimes_{\mathbf R} \widehat{ \mathrm{Sym} } \text{ } \mathrm{H}^2(\mathrm{dP}_n^d,{\mathbf R})^\vee $$
for $n=2,3$. The pairing on $\mathrm{H}^2(\mathrm{dP}_2^d,{\mathbf R})$ induced by the cup product allows the introduction of several familiar differential operators for $n=2$. The ones we will be concerned with are the directional derivative $\nabla_v$, defined for all $v \in \mathrm{H}^2(\mathrm{dP}_2^d,{\mathbf R})$, and the d'Alembert operator $\Box$ associated with the Minkowski metric. We leave the precise definitions in the present circumstances to the reader. 

Returning to the generating function for the number of rational curves through suitable number of points, let ${u}$ be a variable in $\mathrm{H}^2(\mathrm{dP}_n^d,{\mathbf R})$ and define
$$ \Phi_{\mathrm{dP}_n^d}(u,t) = \sum_{\gamma \in \mathrm{H}_2(\mathrm{dP}_n^d,{\mathbf Z})} \left\langle [\mathrm{pt}]^{\deg \gamma - \epsilon_n} \right\rangle_{0,\gamma}^{\mathrm{dP}_n^d} e^{\langle \gamma, u \rangle} t^{\deg{\gamma} } \in F_n^d[[t]],  \eqno(1.3) $$
where $\epsilon_n = 1$ if $n=2$, respectively $0$ if $n=3$. Then (1.2) can be expressed equivalently as $ \smash{ \Phi_{\text{dP}_3^d}\left( u,t \right) = {\mathfrak L} \Phi_{\text{dP}_2^d}\left((j_{2,3}^d)^*u,t \right) }$, where
$$ {\mathfrak L} = \frac{\nabla_\omega \nabla_\omega - d\Box}{d(9-d)}.  \eqno(1.4) $$ 
Here, $\omega = - c_1({\text{dP}_2^d})$, regarded as an element of $\text{H}^2(\text{dP}_2^d,{\mathbf R})$, and $\nabla_\omega \nabla_\omega$ of course stands for the second directional derivative along $\omega$.

\subsection{Acknowledgements} I would like to thank Erwan Brugall\'{e}, Yaim Cooper, Philip Engel, Penka Georgieva, J\'{a}nos Koll\'{a}r, Joe Harris, Quoc Ho, Anand Patel and Alex Perry for useful discussions. One chapter of the author's doctoral thesis is based on this work.

\section{Preliminaries}

\subsection{Genus zero maps to log K3 pairs} In this section and the next, we prove some easy preliminary results which are more convenient to state independently of the rest of the argument. Let $\Sigma = \text{dP}_2^d$ be a complex del Pezzo surface and $E \subset \Sigma$ a smooth anticanonical divisor. By adjunction, $E$ has genus one. We will show, roughly, that there exist only finitely many degree $e$ rational curves on $\Sigma$ which cut $E$ in any predetermined collection of $e$ distinct points. For a homology class $\beta \in \text{H}_2(\Sigma,{\mathbf Z})$, let $\overline{\mathcal M}_{0,e}(\Sigma,\beta)$ be the space of genus zero stable maps to $\Sigma$ with $e$ marked points, where we will choose $e = (E \cdot \beta)$. The evaluation map is $\text{ev}:\overline{\mathcal M}_{0,e}(\Sigma,\beta) \to \Sigma^e$. Let $\xi_1,...,\xi_e \in E$ be distinct closed points on $E$. 

\begin{lem}\label{2.1} The family of curves over $\overline{\mathcal M}_{0,e}(\Sigma,\beta) \times_{\Sigma^e} \{(\xi_1,...,\xi_e)\}$ sweeps out a locus of dimension (at most) one on $\Sigma$. 

If, moreover, for every \textit{proper} subset $I \subset [e] = \{1,2,...,e\}$, we have
$$ {\mathscr O}_E\left(\sum_{i \in I} \xi_i\right) \notin \mathrm{Im}(\mathrm{Pic}(\Sigma) \to \mathrm{Pic}(E)), $$
then $\overline{\mathcal M}_{0,e}(\Sigma,\beta) \times_{\Sigma^e} \{(\xi_1,...,\xi_e)\}$ has dimension zero if nonempty. Moreover, the sources of all the stable maps it parametrizes are smooth.
\end{lem}

\begin{proof} If this was not the case, by properness, this family of curves would have to sweep out the whole surface $\Sigma$. Let $\xi'$ be a point on $E$, different from all the $\xi_i$. Let $(C,f,p_1,...,p_e) \in \overline{\mathcal M}_{0,e}(\Sigma,\beta) \times_{\Sigma^e} \{(\xi_1,...,\xi_e)\}({\mathbf C})$ whose image hits the point $\xi'$. The section $1 \in \text{H}^0(C,f^*{\mathscr O}_\Sigma(E))$ vanishes on at least $e+1$ points in different fibers of $f$, namely $p_1,...,p_e$ and another point $p'$ mapping to $\xi'$, so, by degree considerations, it must vanish on some irreducible component of $C$ which is not contracted by $f$. Hence $C$ has an irreducible component $C_0$ which maps nonconstantly to $E$, which is impossible since $C_0$ has genus zero. 

For the second part, note that the condition implies that for any $(C,f,p_1,...,p_e) \in \overline{\mathcal M}_{0,e}(\Sigma,\beta) \times_{\Sigma^e} \{(\xi_1,...,\xi_e)\}({\mathbf C})$, there is at most one irreducible component of $C$ which is not contracted by $f$. (Here we are implicitly using that $-K_\Sigma = {\mathscr O}_\Sigma(E)$ is ample, so $E$ intersects any divisor on $\Sigma$.) Since $C$ has arithmetic genus zero, it is not hard to deduce that $C$ has no contracted components, so it must in fact be irreducible. Indeed, any contracted component had at most one marked point, so the corresponding vertex in the dual graph has degree at least $2$, but the one-vertex edgeless graph is the unique tree with only one vertex of degree less than $2$ thus proving the claim. Moreover, since $\xi_1,...,\xi_e$ are distinct, $f$ cannot factor as a multiple cover. Then finiteness follows easily from the previous result.
\end{proof}

\begin{rmk} Without further hypotheses, $\overline{\mathcal M}_{0,e}(\Sigma,\beta) \times_{\Sigma^e} \{(\xi_1,...,\xi_e)\}$ may fail to be reduced, e.g. when it contains ramified stable maps, such as the normalization of a cuspidal curve. 
\end{rmk}

\noindent Note that in order for $\overline{\mathcal M}_{0,e}(\Sigma,\beta) \times_{\Sigma^e} \{(\xi_1,...,\xi_e)\}$ to be nonempty, it is necessary that the unique line bundle on $\Sigma$ whose first Chern class is Poincar\'{e} dual to $\beta$ restricts to ${\mathscr O}_E(\xi_1+...+\xi_e)$ on $E$. If that is the case, then the forgetful map 
$$ \overline{\mathcal M}_{0,e}(\Sigma,\beta) \times_{\Sigma^e} \{(\xi_1,...,\xi_e)\} \longrightarrow \overline{\mathcal M}_{0,e-1}(\Sigma,\beta) \times_{\Sigma^{e-1}} \{(\xi_1,...,\xi_{e-1})\} $$ 
is an isomorphism, since the image of any map belonging to the second space necessarily intersects $E$ transversally at $\xi_e$. If, additionally the hypothesis in the second part of lemma \ref{2.1} is satisfied, then
$$ \deg \left[\overline{\mathcal M}_{0,e}(\Sigma,\beta) \times_{\Sigma^e} \{(\xi_1,...,\xi_e)\} \right] = \left\langle[\text{pt}]^{e-1} \right\rangle^{\Sigma}_{0,\beta}. $$
This will be used later in the proof of Theorem 1.2.

\subsection{Lines on Fano threefolds of index two} Let $X = \text{dP}_3^d$ be a smooth Fano threefold of index $2$, Picard rank one and base point free polarization, i.e. $d \in \{2,3,4,5\}$. We will need to know the number $\lambda_{3,1}^d$ of lines (i.e. curves of degree $1$ relative to the polarization by one-half the anti-canonical class) through a general point $x \in X$.

There are several ways to carry out this count, here we sketch a combinatorial approach. Consider the morphism $X \to {\mathbb P}^{d+1}$. This is an embedding if $d \geq 3$, respectively a ramified covering if $d=2$. A general hyperplane section of $X$ is a del Pezzo surface of degree $d$ and it can be described as the blowup of ${\mathbb P}^2$ at $9-d \leq 7$ points, no three of which are collinear and no five of which lie on the same conic. If, $d \geq 3$ let $H \subset {\mathbb P}^{d+1}$ be a general hyperplane containing the projective tangent space to $X$. If $d=2$, we assume that $x$ lies in the ramification locus of $X \to {\mathbb P}^{3}$ and choose $H$ to be the plane tangent to the branch locus at the image of $x$. The case $x \in X$ general will still follow by a semicontinuity argument, which we skip. The pullback $X \times_{{\mathbb P}^{d+1}} H$ has a simple double point at $x$ and is smooth elsewhere. Moreover, it will contain any line on $X$ passing through $x$. Resolving this singularity, we obtain a surface $S$ with a $(-2)$-curve. The lines on $X$ through $x$ correspond to $(-1)$-curves on $S$ intersecting the $(-2)$-curve.

The surface $S$ also can be described as the blowup of ${\mathbb P}^2$ at $9-d$ points, but with the feature (for instance) that $3$ of the points have become collinear. Of course, the $(-2)$-curve is simply the proper transform of the line through these $3$ points. Let $A,B,C \in {\mathbb P}^2$ be the three collinear points and $P_1,...,P_{6-d} \in {\mathbb P}^2$ the remaining $6-d$ points. The only $(-1)$-curves intersecting the $(-2)$-curve are the three exceptional divisors of the blow up at $A$, $B$ and $C$, the proper transforms of the lines $P_iP_j$, $i<j$, and the proper transforms of conics passing through one of the points $A,B,C$ and $4$ of the $P_i$, so the answer is
$$ \lambda_{3,1}^d = 3{6-d \choose 4} + {6-d \choose 2} + 3 $$
which is $12$ for $d=2$, $6$ for $d=3$, $4$ for $d=4$ respectively $3$ for $d=5$. Note, in particular, that these counts are compatible with the formula in theorem 1.2.

\subsection{Commutative algebra preliminaries} In this subsection, we prove a standard commutative algebra result, which will be used to formalize in families the trick involving Bezout's theorem, which was outlined in the introduction. The reader may skip ahead if he or she wishes and return to this section when its purpose becomes clear. 

\begin{lem}\label{2.3} Let $f:(A,{\mathfrak m}) \to (B, {\mathfrak n})$ be a flat local homomorphism of local Noetherian rings. If ${\mathfrak I} \subseteq {\mathfrak m}^e$ is an ideal of $B$ such that $B/{\mathfrak I}$ is flat as an $A$-module, then ${\mathfrak I} = (0)$.
\end{lem}

\begin{proof} First, tensoring the short exact sequence $0 \to {\mathfrak m} \to A \to A/{\mathfrak m} \to 0$ with $B$ we obtain an exact sequence of $A$-modules
$$ 0 \longrightarrow {\mathfrak m} \otimes_A B \longrightarrow A \otimes_A B \longrightarrow A/{\mathfrak m} \otimes_A B \longrightarrow 0. $$
The middle term is $B$ and the next term is $B/{\mathfrak m}B = B/{\mathfrak m}^e$, so flatness of $B$ over $A$ implies that the $A$-module homomorphism $B \otimes_A {\mathfrak m} \to {\mathfrak m}^e$ is actually an isomorphism. Second, tensoring the short exact sequence $0 \to {\mathfrak m} \to A \to A/{\mathfrak m} \to 0$ with $B/{\mathfrak I}$ we obtain an exact sequence of $A$-modules
$$ 0 \longrightarrow {\mathfrak m} \otimes_A B/{\mathfrak I} \longrightarrow A \otimes_A B/{\mathfrak I} \longrightarrow A/{\mathfrak m} \otimes_A B/{\mathfrak I} \longrightarrow 0. $$
First, $A/{\mathfrak m} \otimes_A B/{\mathfrak I}$ is $(B/{\mathfrak I})/{\mathfrak m}(B/{\mathfrak I}) = (B/{\mathfrak I})/{\mathfrak m}^e(B/{\mathfrak I})$. The kernel of the surjective composition $B \to B/{\mathfrak I} \to (B/{\mathfrak I})/{\mathfrak m}^e(B/{\mathfrak I})$ is just ${\mathfrak m}^e$ hence $A/{\mathfrak m} \otimes_A B/{\mathfrak I} \cong B/{\mathfrak m}^e$. We therefore have an isomorphism ${\mathfrak m} \otimes_A B/{\mathfrak I} \cong {\mathfrak m}^e/{\mathfrak I}$. 

However, we have a well defined surjective $A$-module homomorphism ${\mathfrak m} \otimes_A B/{\mathfrak I} \to {\mathfrak m}^e/{\mathfrak m}^e{\mathfrak I}$, which fits in the following commutative diagram of $A$-modules. 
\begin{center}
\begin{tikzpicture}

\matrix [column sep  = 17mm, row sep = 15mm] {
	\node (nw) {${\mathfrak m}^e$}; & 
	\node (nc) {${\mathfrak m} \otimes_A B$}; &
	\node (ne) {${\mathfrak m}^e$}; \\
	\node (sw) {${\mathfrak m}^e/{\mathfrak m}^e{\mathfrak I}$}; & 
	\node (sc) {${\mathfrak m} \otimes_A (B/{\mathfrak I})$}; &
	\node (se) {${\mathfrak m}^e/{\mathfrak I}$}; \\
};

\draw[-{Stealth[length=2mm]}, thin] (nc) -- (nw);
\draw[-{Stealth[length=2mm]}, thin] (nc) -- (ne);
\draw[-{Stealth[length=2mm]}, thin] (sc) -- (sw);
\draw[-{Stealth[length=2mm]}, thin] (sc) -- (se);
\draw[-{Stealth[length=2mm]}, thin] (nw) -- (sw);
\draw[-{Stealth[length=2mm]}, thin] (nc) -- (sc);
\draw[-{Stealth[length=2mm]}, thin] (ne) -- (se);

\end{tikzpicture}
\end{center}
The vertical lateral maps are simply quotients. The lower left horizontal map is surjective and all the other horizontal maps are isomorphisms. It follows that ${\mathfrak I}={\mathfrak m}^e{\mathfrak I}$. Inductively, ${\mathfrak I}={({\mathfrak m}^e)}^k{\mathfrak I}$ for all positive integers $k$, hence
$$ {\mathfrak I} \subseteq \bigcap_{k \geq 0} {({\mathfrak m}^e)}^k  = (0), $$
by Krull's intersection theorem, so ${\mathfrak I}=(0)$, as desired. \end{proof}

\begin{cor}\label{2.4} As above, let $f:(A,{\mathfrak m}) \to (B, {\mathfrak n})$ be a flat local homomorphism of local Noetherian rings. Let $M$ be a $B$-module such that:

$\bullet$ $M$ is flat as an $A$-module; and

$\bullet$ $B/{\mathfrak m}^e \otimes_B M$ is a free rank one module over $B/{\mathfrak m}^e$.

\noindent Then $M$ is a free rank one module over $B$.
\end{cor}

\begin{proof} 
This follows immediately from Lemma \ref{2.3} and Nakayama's lemma. Indeed, by Nakayama's lemma, $M$ has one generator over $B$, i.e. $M \cong B/{\mathfrak I}$ for some ideal ${\mathfrak I} \subseteq B$. Then the second bullet reads $B/{\mathfrak m}^e \otimes_B B/{\mathfrak I} = B/({\mathfrak m}^e + {\mathfrak I})$ is a free rank one $B/{\mathfrak m}^e$-module, hence ${\mathfrak I} \subseteq {\mathfrak m}^e$. Thus, we may invoke \ref{2.3}. 
\end{proof}

\section{Homology classes on surfaces in a pencil}

\subsection{The Lefschetz pencil} As in section 2.2, we consider a smooth Fano threefold $X$ of index $2$ and degree $d$, with polarization ${\mathscr O}_X(1)$. Assume that $\text{Pic}(X) = {\mathbf Z}$ and the polarization is base-point free, which amounts to $d \in \{2,3,4,5\}$. Consider a general pencil of sections of the polarizing line bundle with base locus $E$ and total space $\rho:W \to {\mathbb P}^1$, where $W$ is the blowup of $X$ along $E$. By adjunction, $p_a(E)=1$. The members of the pencil are generically smooth del Pezzo surfaces of the same degree $d$ as $X$. Let ${\mathbb P}^\circ \subset {\mathbb P}^1$ parametrize smooth del Pezzo fibers and $W^\circ \subset W$ its preimage. Choose a closed point $b \in {\mathbb P}^\circ$ and $\beta_b \in \text{H}_2(W_b,{\mathbf Z})$. The Poincar\'{e} dual of $\beta_b$ is the first Chern class of a (uniquely determined) line bundle ${\mathscr L}_b$ on the surface $W_b$. 

First, we introduce some notation. We will consider objects $(N,\langle \cdot,\cdot \rangle, \nu)$ consisting of the following data: (1) a finitely generated free abelian group $N$; (2) a bilinear map $\langle \cdot,\cdot \rangle: N \times N \to {\mathbf Z}$; and (3) a nonzero distinguished element $\nu \in N$. A morphism between two such objects $(N,\langle \cdot,\cdot \rangle_1, \nu)$ and $(M,\langle \cdot,\cdot \rangle_2, \mu)$ is a map $\varphi:N \to M$ such that $\varphi(\nu) = \mu$ and $\langle \varphi(v),\varphi(w) \rangle_2 = \langle v,w \rangle_1$. We denote the resulting category by ${\mathfrak D}$. The purpose of this category is merely to simplify language.

Set $r=9-d$. Let $(H_r,\langle \cdot,\cdot \rangle, \omega)$ be an object of ${\mathfrak D}$ defined by $H_r = {\mathbf Z}\ell_0 \oplus {\mathbf Z}\ell_1 \oplus ... \oplus {\mathbf Z}\ell_r$,
$$ \langle \ell_i, \ell_j \rangle=
\begin{cases}
\hphantom{-} 0  & \text{if $i \neq j$},\\
\hphantom{-} 1 &\text{if $i=j=0$},\\
-1 &\text{if $i=j \geq 1$}
\end{cases} \eqno(3.1)
$$ 
and $\omega = -3\ell_0+\ell_1+\ell_2+...+\ell_r$. Note that this is isomorphic to $\text{H}^2(\Sigma,{\mathbf Z})$ of a degree $d$ del Pezzo surface with the cup product and canonical class. The abelian subgroup $\omega^\perp = \{v: \langle v,\omega\rangle = 0 \}$ with the negative pairing is the lattice $E_r$, if $d \in \{2,3,4,5\}$. We are abusing notation by writing $E_r$ for $r =4, 5$ instead of $A_4$ and $D_5$ respectively. The complexification $E_r \otimes {\mathbf C} = \omega^\perp_{\mathbf C}$ is the (dual) Cartan subalgebra ${\mathfrak h}^\vee$ of the corresponding Lie algebra and the restriction of the pairing is the Killing form. Of course, the Killing form identifies ${\mathfrak h}^\vee$ and ${\mathfrak h}$, so we will simply write ${\mathfrak h}$ despite the fact that we sometimes mean the dual. Recall that the automorphism group $\text{Aut}_{\mathfrak D}(H_r,\langle \cdot,\cdot \rangle, \omega)$ is the Weyl group $G={\mathcal W}(E_r)$. Another folklore fact which will be used is the following: up to scalars, the Killing form is the only $G$-invariant bilinear map on ${\mathfrak h}$.

Consider the monodromy action $\smash{\pi_1({\mathbb P}^\circ,b) \to \text{Aut}_{\mathfrak D}(\mathrm{H}^2(W_b,{\mathbf Z}),\cup,c_1({\mathscr T}^\vee_{W_b}))}$. As $W \to {\mathbb P}^1$ is a Lefschetz pencil, it is well-known that the image of the monodromy map is the full Weyl group $G$. Let $B^\circ$ parametrize pairs $(t,\varphi_t)$ consisting of a closed point $t \in {\mathbb P}^\circ$ and a ${\mathfrak D}$-isomorphism of $(\text{H}^2(W_t,{\mathbf Z}),\cup,c_1({\mathscr T}^\vee_{W_t}))$ with $\smash{(H_r,\langle \cdot , \cdot \rangle,\omega)}$. Then $B^\circ \to {\mathbb P}^\circ$ is finite and \'{e}tale. A sketch of an algebraic construction of $B^\circ$ goes as follows. Let $F \to {\mathbb P}^\circ$ be the relative Fano scheme of lines in the smooth fibers of the pencil. The locus 
$$ \{ (\ell_1,...,\ell_r): \ell_i \cap \ell_j = \O, \text{ for }i \neq j \} \subset F \times_{{\mathbb P}^\circ} F \times_{{\mathbb P}^\circ} ... \times_{{\mathbb P}^\circ} F $$
parametrizing $r$-tuples of mutually disjoint lines in the same fiber is both open and closed. We define $B^\circ$ to be this locus and $B$ to be its nonsingular completion. By properness, the \'{e}tale map $B^\circ \to {\mathbb P}^\circ$ extends to a (ramified) finite map $B \to {\mathbb P}^1$.

Let $W^\circ_B$ be the pullback of $W$ to $B^\circ$ and $W_B$ the pullback of $W$ to $B$. Although we aren't going to use it, we should point out that all singular points of $W_B$ are ordinary double points so, after small resolutions, the family $W_B \to B$ becomes a smooth family $Y \to B$. Note, nevertheless, that this family would still be not a topological fibration, since the new fibers have other intersection forms. By construction, the relative Picard scheme $\text{Pic}(W^\circ_B/B^\circ)$ is simply $H_r \times B^\circ$, so every $\beta \in H_r$ naturally induces a section $\sigma(\beta):B^\circ \to \text{Pic}(W^\circ_B/B^\circ)$. Before stating the main enumerative problem of this section, we prove a lemma which will turn out very useful later in ruling out potential multiplicities.

\begin{lem}\label{3.1} The composition $B^\circ \to \mathrm{Pic}(W^\circ_B/B^\circ) \to \mathrm{Pic}(E)$ of $\sigma(\beta)$ with the restriction map is constant if and only if $\beta$ is a rational multiple of $\omega$.
\end{lem}

\begin{proof} The if direction is trivial since $-K_{W_t} = {\mathscr O}_{W_t}(1)$ restricts to ${\mathscr O}_E(1)$ on $E \subset W_t$. To prove the converse, we use monodromy. Let $\Lambda \subset H_r$ be the set of all $\lambda \in H_r$ for which the composition in the statement of the lemma is actually constant. Of course, $\Lambda \neq \O$ since it contains $\omega$. It is not hard to check that $\Lambda$ has the following properties: (1) $\Lambda$ is a subgroup of $H_r$, (2) $\Lambda$ is $G$-invariant due to monodromy and (3) if $m \neq 0$ is an integer and $m \lambda \in \Lambda$, then $\lambda \in \Lambda$.

An easy exercise in lattice theory proves that the properties above imply that either $\Lambda$ consists precisely of the multiples of $\omega$, or $\Lambda = H_r$. We claim that the latter is impossible. Indeed, we will show that $\ell=\ell_1 \notin \Lambda$. Geometrically, $\ell$ is the class of a line on a del Pezzo surface of degree $d$. In section {2.2}, we've shown that there exist lines passing through each point of the original threefold $X$. Applying this to points of $E$ and noting that any line on $X$ intersecting $E$ has to lie inside some $W_t$, we conclude that $\ell \notin \Lambda$. \end{proof}

\noindent We will also use the infinitesimal version of the lemma: the differential of the composition above is generically nonzero, unless $\beta \in {\mathbf Q} \omega$. 

\subsection{Counting homology classes} Let $\beta \in H_r$ be an element corresponding to the chosen $\beta_b \in \text{H}_2(W_b,{\mathbf Z})$ under a suitable isomorphism. Fix ${\mathscr L}_E \in \text{Pic}^e(E)$. In this section, we want to address the following enumerative problem. 

\begin{prob}\label{3.2} Assuming sufficiently general choices, how many pairs $(t,{\mathscr L}_t)$ consisting of a closed point $t \in {\mathbb P}^\circ$ and ${\mathscr L}_t \in \mathrm{Pic}(W_t) \cong \mathrm{H}^2(W_t,{\mathbf Z})$ such that

$\bullet$ ${\mathscr L}_t$ restricts on $E$ to the line bundle ${\mathscr L}_E$; and

$\bullet$ there is a $\mathfrak D$-isomorphism $(\mathrm{H}^2(W_t,{\mathbf Z}),\cup,c_1({\mathscr T}_{W_t}^\vee)) \cong (H_r,\langle \cdot , \cdot \rangle,\omega)$ mapping ${\mathscr L}_t \mapsto \beta$ are there? \end{prob}

\noindent The way we will answer 3.2 is similar to the way most enumerative questions are answered: by doing intersection theory on a moduli space. Consider the functor $\text{Sch}_{\mathbf C}^\text{op} \to \text{Set}$ mapping a scheme $S$ over ${\mathbf C}$ to the set of homogeneous $S$-group scheme homomorphisms  $H_r \times  S \to \text{Pic}(E) \times S$ which send $\omega$ to ${\mathscr O}_E(-1)$ fiberwise. We require this to be homogeneous with respect to the grading on $H_r$ given by pairing with $-\omega$ and the natural grading by degree on $\text{Pic}(E)$. The map on arrows is defined in the obvious way by pullback.

It is not hard to prove that the functor defined above is represented by an abelian variety $A$ of complex dimension $r$. The tangent bundle of $A$ is naturally isomorphic to ${\mathfrak h} \otimes_{\mathbf C} {\mathscr O}_A = E_r \otimes_{\mathbf Z} {\mathscr O}_A$.  If we defined a similar functor without the requirement that $\omega \mapsto {\mathscr O}_E(-1)$ but not dropping the homogeneity condition, we'd have obtained an abelian variety of dimension $r+1$ which we'll denote by $V$. It is clear by definition that $A$ sits naturally inside $V$. 

Recall that if $T$ is a complex abelian variety, or more generally a complex analytic torus, then the isomorphism of graded-commutative ${\mathbf C}$-algebras
$$ \left( \mathrm{H}^*(T,{\mathbf C}), \cup \right) \cong \left({\bigwedge} \mathrm{H}^1(T,{\mathbf C}), \wedge \right) \eqno(3.2) $$
respects the Hodge structures on both sides. In particular, there is a canonical isomorphism $\text{H}^{1,1}(T) \cong \text{H}^{1,0}(T) \otimes \text{H}^{0,1}(T)$.

The necessary moduli space for solving Problem \ref{3.2} directly is $[A/G]$, but to avoid being unnecessarily sophisticated, we will work simply with $A$.  The family $E \times B^\circ \subset W_B^\circ \to B^\circ$ induces a morphism $B^\circ \to A$. By properness, this map is canonically extended to a $G$-equivariant morphism $\tilde{\alpha}:B \to A$, which fits in a commutative square with the analogous map $\alpha: {\mathbb P}^1 \to [A/G]$ and the quotients by $G$. 

The evaluation of a homogeneous homomorphism $H_r \to \text{Pic}(E)$ at $\beta$ is a morphism of abelian varieties $\hat{\beta}:A \to \text{Pic}^e(E)$. An equivalent formulation of Problem \ref{3.2} is to ask what is the number of triples $(t,{\mathscr L}_t,\varphi_t)$ where $\varphi_t$ is an isomorphism as in the second bullet and ${\mathscr L}_t|_E \cong {\mathscr L}_E$. The theoretical answer to this problem is $\smash{ \deg \hat{\beta} \tilde{\alpha} }$. To obtain the answer in the original formulation, simply divide by $|\text{Stab}_G(\beta)|$.

Consider the pullback on cohomology $\smash{ \tilde{\alpha}^*:\text{H}^{1,1}(A) \to \text{H}^{1,1}(B) }$. The domain of definition is canonically isomorphic to ${\mathfrak h} \otimes \overline{\mathfrak h}$, so we obtain a $G$-invariant map 
$$ \tilde{\alpha}^*:{\mathfrak h} \otimes \overline{\mathfrak h} \to \text{H}^{1,1}(B). $$ 
The crucial observation is that by the uniqueness up to scalars of $G$-invariant bilinear maps, $\smash{ \tilde{\alpha}^* }$ factors through the map $v \otimes \overline{w} \mapsto \langle v, w \rangle$, where the inner product is the Killing form. Note that we are actually using the general fact stated above for sesquilinear maps, which is equivalent since the action of $G$ on ${\mathfrak h}$ is real, i.e. it preserves the ${\mathbf R}$-span of $E_r$. It follows that
$$ \tilde{\alpha}^*(v \otimes \overline{w}) = - \langle v, w \rangle \tilde{\alpha}^* (\ell_1 \otimes \ell_1), \eqno(3.3) $$
so we've boiled down the calculation to computing $\tilde{\alpha}^* (\ell_1 \otimes \ell_1) \in \text{H}^{1,1}(B)$ and $\smash{\hat{\beta}^*([\text{pt}])}$.

\begin{lem}\label{3.3} Let $\beta^\perp $ be the projection of $\beta$ to the orthogonal complement of $\omega$. The pullback morphism $\mathrm{H}^{1,1}(\mathrm{Pic}^e(E)) \to \mathrm{H}^{1,1}(A)$ maps the Poincar\'{e} dual of the class of a point to the pure tensor $\beta^\perp \otimes \beta^\perp \in {\mathfrak h} \otimes \overline{\mathfrak h} \cong \mathrm{H}^{1,1}(A)$. \end{lem}

\begin{proof} First, we should point out that all $\text{H}^{i,j}(\text{Pic}^e(E))$ can be naturally identified with ${\mathbf C}$, since $\text{H}^0(E, {\mathscr O}_E) \cong {\mathbf C}$. Consider the pullback maps
$$ \hat{\beta}_{i,j}^*:\text{H}^{i,j}(\text{Pic}^e(E)) \longrightarrow \text{H}^{i,j}(A) $$
for all $i,j \in \{0,1\}$. Note that $\smash{ \hat{\beta}_{1,1}^* = \hat{\beta}_{1,0}^* \otimes \hat{\beta}_{0,1}^* }$ and the two factors are naturally complex conjugate. Consider the dual differential 
$$ d\hat{\beta}^\vee: \hat{\beta}^*\Omega^1\text{Pic}^e(E) \longrightarrow \Omega^1(A). $$
This map is in fact simply the morphism ${\mathscr O}_A \to {\mathscr O}_A \otimes {\mathfrak h}$ given by $s \mapsto s \otimes \beta^\perp$. This will be proved shortly, let's just assume it for now. Taking global sections of $d\hat{\beta}^\vee$ we recover the map $\hat{\beta}_{1,0}^*$, which was therefore the unique linear map ${\mathbf C} \to {\mathfrak h}$ such that $1 \mapsto \beta^\perp$. Having understood $\smash { \hat{\beta}_{1,0}^* }$, we extrapolate to $\smash {\hat{\beta}_{0,1}^*}$ by conjugation and then to $\smash {\hat{\beta}_{1,1}^*}$ by tensoring, obtaning the desired conclusion. \end{proof}

\begin{claim}\label{3.4} The map $d\hat{\beta}^\vee$ is just ${\mathscr O}_A \to {\mathscr O}_A \otimes {\mathfrak h}$, $s \mapsto s \otimes \beta^\perp$. \end{claim}

\begin{proof} Note that there is a related map $d\hat{\beta}_V^\vee: \hat{\beta}_V^*\Omega^1\mathrm{Pic}^e(E) \to \Omega^1(V)$. This is easily seen to be simply tensoring with $\beta$. To obtain the original $d\hat{\beta}^\vee$, we need to restrict to $A$ and compose this map with the restriction $\Omega^1(V) \otimes {\mathscr O}_A \to \Omega^1(A)$, which is simply the sheafified version of the projection map $H_r \otimes {\mathbf C} \to {\mathfrak h}$ to the orthogonal complement of $\omega$. \end{proof}

\noindent Now we can solve 3.2. Denote $\Delta(\omega,\beta) = (\omega \cdot \beta)^2 - \left(\omega^2\right)\left(\beta^2\right)$. Note that $\Delta(\omega,\beta) \geq 0$ by the Hodge index theorem, with equality if and only if $\beta$ is a multiple of $\omega$. By (3.3) and Lemma \ref{3.3}, we get
$$ \deg \hat{\beta} \tilde{\alpha} = \tilde{\alpha}^*\left( \hat{\beta}^*\left([\text{pt}] \right) \right) = \tilde{\alpha}^*\left( \beta^\perp \otimes \beta^\perp \right) = - \frac{\Delta(\omega,\beta)}{(\omega^2)} \tilde{\alpha}^*\left( \ell_1 \otimes \ell_1 \right), \eqno(3.4) $$
for all $\beta \in H_r$. Plugging $\beta = \ell_1$ in $(3.4)$, we obtain $\smash{ \deg ( \hat{\ell}_1 \tilde{\alpha} ) = -(1+1/d)\tilde{\alpha}^*\left( \ell_1 \otimes \ell_1 \right)}$. However, recall that all lines on $X$ intersecting $E$ actually lie in one of the fibers $W_t$ of the pencil. Therefore, the argument so far for $e=1$ shows that
$$ \lambda^d_{2,0} \deg \hat{\ell}_1 \tilde{\alpha} =  |G| \lambda^d_{3,1}, \text{ hence } \tilde{\alpha}^*\left( \ell_1 \otimes \ell_1 \right) = -|G|\frac{d\lambda_{3,1}^d}{(d+1)\lambda_{2,0}^d}, \eqno(3.5) $$
where $\lambda^d_{2,0} = |G\cdot \ell_1|$ is the number of lines on a degree $d$ del Pezzo surface and $\lambda^d_{3,1}$ is, as in section {2.2}, the number of lines through a point on a degree $d$ del Pezzo threefold. 
\begin{center}
\begin{tabular}{l|lllll}
Type & $d$ & $\lambda^d_{2,0}$ & $\lambda^d_{3,1}$ & $|G|$ & $\frac{d\lambda_{3,1}^d}{(d+1)\lambda_{2,0}^d}$ \\ \hline
double cover ramified along a quartic & $2$ & $56$ & $12$ & $2903040$ & $1/7$ \\
cubic threefold & $3$ & $27$ & $6$ & $51840$ & $1/6$ \\
$(2,2)$-complete intersection & $4$ & $16$ & $4$ & $1920$ & $1/5$ \\
plane section of Grassmannian ${\mathbb G}(1,4)$ & $5$ & $10$ & $3$ & $120$ & $1/4$ \\
\end{tabular}
\end{center}
The final entry is $1/r$ in all four cases. Therefore, by (3.4) and (3.5), we get
$$ \deg \hat{\beta} \tilde{\alpha} = |G|\frac{\Delta(\omega,\beta)}{d(9-d)}, $$
from which we obtain the answer by correcting with the stabilizer factor. 

Finally, we will rearrange the enumerative problem we've just solved in a more convenient form. Let $S({\mathscr L}_E,\beta)$ be the set of solutions to question \ref{3.2}. As a scheme, this lives on $\text{Pic}(W^\circ/{\mathbb P}^\circ)$ for general ${\mathscr L}_E$, since there are only countably many divisor classes on $E$ obtained by restricting divisor classes on the singular fibers of $W \to {\mathbb P}^1$. Moreover, it is reduced by Lemma \ref{3.1} and its infinitesimal version, since no ${\mathscr L}_t$ can be a multiple of $K_{W_t}$ if ${\mathscr L}_E(-m)$ is not torsion for all integers $m$. In the curve counting problem, we will encounter a $G$-constant counting function $N:H_r \to {\mathbf N}$. (In that case, $N$ counts the number of class $\beta$ genus $0$ stable maps through $e-1 = \deg \beta -1$ points on a degree $d$ del Pezzo surface $\Sigma$.) Then, as purely combinatorial statement, we have
$$ \sum_{G\beta \in H_r^e/G} |S({\mathscr L}_E,\beta)| {N(\beta)} = \sum_{\beta \in H^e_r} \frac{|G|\Delta(\omega,\beta)}{d(9-d)} \frac{1}{|\text{Stab}_G(\beta)|} \frac{N(\beta)}{|G\cdot \beta|} = $$
$$ = \frac{1}{d(9-d)}\sum_{\beta \in H^e_r}\Delta(\omega,\beta)N(\beta), \eqno(3.6) $$
where $H_r^e$ is the degree $e$ piece of $H_r$.

\section{Interpolating points on a genus one curve}

Let $X$ be a del Pezzo variety of dimension $n$ and degree $d \in \{2,3,4,5,8\}$, a homology class $\beta \in \text{H}_2(X,{\mathbf Z}) $ and $\overline{\mathcal M}_{0,0}(X,\beta)$ the space of class $\beta$ genus $0$ stable maps to $X$. It is easy to check that that $\text{H}_2(X,{\mathbf Z})$ is torsion free. The hyperplane class ${\mathscr O}(1) \in \mathrm{Pic}(X)$ is uniquely determined by the property $-K_X = {\mathscr O}(n-1)$. Let $e = (\beta \cdot {\mathscr O}(1))$ be the degree of $\beta$. 

Recall from the introduction that we want to analyze the incidence correspondence between $m$-tuples of points on $X$ and rational curves containing them. In the stable map compactification, the correspondence translates naturally as the evaluation map $\text{ev}:\overline{\mathcal M}_{0,m}(X,\beta) \to X^m$. Roughly, the evaluation map is dominant if and only if it is possible to interpolate $m$ general points on $X$ with a rational curve of class $\beta$. This is expected to happen when $\text{vdim }\overline{\mathcal M}_{0,m}(X,\beta) \geq \dim X^m$, or equivalently,
$$ m \leq \frac{(-K_X \cdot \beta) - 2}{\dim X - 1} + 1 = \frac{(n-1)e - 2}{n-1} + 1, \eqno(4.1) $$
since $\text{vdim }\overline{\mathcal M}_{0,m}(X,\beta) = \dim X -(K_X \cdot \beta) + m - 3$, so the largest integer $m = m_{\text{max}}$ for which the inequality $(4.1)$ holds is
$$ m_{\text{max}}=
\begin{cases}
\hphantom{-} e  & \text{if $n \geq 3$},\\[2ex]
e-1 &\text{if $n = 2$}.
\end{cases} \eqno(4.2)
$$ 
Before outlining the approach, we introduce some notation. 

Let $X$ be any smooth variety, $\beta$ a homology class representing a curve class, $Y$ a closed subvariety of $X$ and $m \leq n$. We write $\smash{ \overline{\mathcal M}^{m}_{g,n} (X,Y;\beta) = \overline{\mathcal M}_{g,n} (X,\beta) \times_{X^m} Y^m }$, where the map from $\smash{ \overline{\mathcal M}_{g,n} (X,\beta)}$ to $X^m$ is $(\text{ev}_1,...,\text{ev}_m)$. If $U \subset Y^m$ is open, let  
$$ \overline{\mathcal M}^{m}_{g,n} (X,Y;\beta)|_U = \overline{\mathcal M}_{g,n} (X,Y;\beta) \times_{Y^m} U. $$
If $m=n$, we drop the superscript. If $n=m+1$, we can think of $\smash{ \overline{\mathcal M}^{m}_{g,m+1} (X,Y;\beta) }$ as the universal curve over $\smash{ \overline{\mathcal M}_{g,m} (X,Y;\beta) }$ and we denote it by $\smash{ \overline{\mathcal C}_{g,m} (X,Y;\beta)}$. If $g=0$, which is the only case treated in this paper, we drop the subscript indicating the genus. 

Let us return to the problem. Let $E \subset X$ be a section of $X$ by $n-1$ general hyperplanes. By adjunction and Bertini, the property of being del Pezzo is preserved at each step, so $E$ is a smooth genus one curve. Set 
$$ V = \mathrm{H}^0(X,{\mathscr I}_{E/X}(1)) \subset \mathrm{H}^0(X,{\mathscr O}(1)). $$ 
Roughly, the main observation is the following: a curve of class $\beta$ which meets $E$ at $e$ distinct points is forced to lie in a surface on $X$ containing $E$, obtained by cutting $X$ with $n-2$ hyperplanes. To use this observation, we have to formalize it in families. Let $S$ be any finite type scheme over ${\mathbf C}$ and
$$ (C,\pi,f,p_1,p_2,...,p_e) \in \overline{\mathcal M}_e(X,E;\beta)|_{\Delta^c}(S), $$
where $\Delta^c = E^e \backslash \Delta$ and $\Delta$ is the big diagonal of $E^e$. Denote by $D_i \subset C$ the image of the closed embedding $p_i:S \to C$ and set $D = \sum D_i$. We start by proving the observation for of stable maps, then proceed with the formalism in families. 

Much of the formalism below is forced by the possibility that $S$ is not reduced. Since we can't say a priori that $\overline{\mathcal M}_e(X,E;\beta)$ is generically reduced (i.e. that it is not contained in the ramification locus of the evaluation map), this is a difficulty we are forced to face. The generic smoothness of $\overline{\mathcal M}_e(X,E;\beta)$ will be a corollary of the subsequent analysis.

\begin{lem}\label{4.1} If $S=\text{Spec}({\mathbf C})$, then $f^*{\mathscr O}(1) \cong {\mathscr O}_C(D)$. Moreover, there is no irreducible component of $C$ mapped constantly to a point on $E$.\end{lem}

\begin{proof} Let $M$ be a finite set indexing all \textit{maximal connected} curves of arithmetic genus zero $C'_\mu \subset C$ which are contracted by $f$ to a point on $E$ and let $C = \bigsqcup_{\mu \in M}C'_\mu \cup C_0$ with $C_0$ possibly disconnected. The dual graph $\Gamma_\mu$ of each $C'_\mu$ is a tree. We further decorate each dual graph with "legs" for each intersection point with $C_0$. For all $\mu \in M$, let $\nu_\mu$ be the number of $i$ such that $p_i \in C'_\mu$ and $\lambda_\mu$ the number of legs in the dual graph $\Gamma_\mu$. Since $f(p_i) \neq f(p_j)$ for $i \neq j$, $\nu_\mu \leq 1$ for all $\mu \in M$. By stability, this implies that no vertex of any $\Gamma_\mu$ can be a leaf of the dual graph of $C$. Therefore, there is at least one leg attached to each leaf of $\Gamma_\mu$, so $\lambda_\mu \geq 2$. In particular, $\lambda_\mu > \nu_\mu$. Similarly, we let $\nu_0$ be the number of $i$ such that $p_i \in C_0$ and $\lambda_0 = \sum \lambda_\mu$. 

Let $D_0$ be the restriction of $D$ to $C_0$ and $D_\lambda$ the divisor on $C_0$ consisting of the $\lambda_0$ "bridge points" to the union of the components on which $f$ maps constantly to $E$. Clearly, $D_0$ and $D_\lambda$ are reduced and have disjoint supports. Let $f_0$ be the restriction of $f$ to $C_0$. We claim that the line bundle 
$$ {\mathscr L}_0:= f_0^*{\mathscr O}(1) \otimes {\mathscr O}_{C_0}(-D_0-D_\lambda) $$
admits sections with finitely many zeroes. We may argue using the pullback map on global sections $\text{H}^0(X,{\mathscr I}_{E/X}(1)) \to \text{H}^0(C_0,{\mathscr L}_0)$. Indeed, we may choose a section of ${\mathscr I}_{E/X}(1)$ which is not identically zero on the image of any component of $C_0$ and correspondingly map it to a section of ${\mathscr L}_0$ with finitely many zeroes. From here, we get the inequality $e - \nu_0 - \lambda_0 \geq 0$, or $\nu_0 +\mu_0 \leq e$. Since $\nu_0 + \sum \nu_\mu = e$, it follows that 
$$ \sum \lambda_\mu = \lambda_0 \leq \sum \nu_\mu, $$
which contradicts $\lambda_\mu > \nu_\mu$ for all $\mu \in M$, unless $M = \varnothing$. Regardless, the line bundle ${\mathscr L}_0$, now $f^*{\mathscr O}(1) \otimes {\mathscr O}_C(-D)$, still has sections with finitely many zeroes and visibly has degree $0$, so it must be trivial. \end{proof}

\begin{lem}\label{4.2} As above, define ${\mathscr L}= f^*{\mathscr O}(1) \otimes {\mathscr O}_C(-D)$. Then $\pi_*{\mathscr L}$ is invertible. Moreover, if $ \varphi_S:V \otimes {\mathscr O}_S \to \pi_*{\mathscr L}$ is the natural ${\mathscr O}_S$-modules map, there exists a unique morphism $\psi_S:S \to {\mathbb P}V$ such that $\mathrm{Ker}(\varphi_S) = \psi_S^* {\mathscr U}_{{\mathbb P}V}$, where ${\mathscr U}_{{\mathbb P}V}$ is the tautological subbundle of $V \otimes {\mathscr O}_{{\mathbb P}V}$.\footnote{We are using the Grothendieck convention for projective spaces, which says that the closed points of ${\mathbb P}V$ correspond to codimension one, rather than dimension one, subspaces of $V$.}\end{lem}

\begin{proof} First, let us spell out the construction of the map $\varphi_S$ in the statement. The ${\mathscr O}_X$-modules homomorphism $\text{H}^0({\mathscr O}(1)) \otimes {\mathscr O}_X \to {\mathscr O}(1)$ pulls back via $f$ to a map $\text{H}^0({\mathscr O}(1)) \otimes {\mathscr O}_C \to f^*{\mathscr O}(1)$. By the adjoint property, we get an ${\mathscr O}_S$-modules homomorphism 
$$ \text{H}^0({\mathscr O}(1)) \otimes {\mathscr O}_S \longrightarrow \pi_*f^*{\mathscr O}(1). $$ 
This map composed further with $\pi_*f^*{\mathscr O}(1) \to \pi_*(f^*{\mathscr O}(1)|_D)$ vanishes on $V \otimes {\mathscr O}_S$, so it induces an ${\mathscr O}_S$-modules homomorphism $\varphi_S:V \otimes {\mathscr O}_S \to \pi_*(f^*{\mathscr O}(1) \otimes {\mathscr O}_C(-D)) = \pi_*{\mathscr L}$. For any closed point $s \in S$, the map 
$$ V \otimes_{\mathbf C} \kappa(s) \longrightarrow \pi_*{\mathscr L} \otimes_{{\mathscr O}_S} \kappa(s) \longrightarrow \text{H}^0(C_s,{\mathscr L}_s) $$ 
is nonzero; otherwise, the image of $f_s$ would be contained inside $E$, which is impossible. However, ${\mathscr L}_s$ is trivial by lemma \ref{4.1} meaning that the rightmost term above is $1$-dimensional, so the composed map above is surjective. 

Furthermore, the second map in the composition has to be surjective as well, so by the cohomology and base change theorem, it is actually an isomorphism. By the same theorem, this property extends automatically to the non-closed points. The $\pi$-pushorward of any torsion-free sheaf on $C$ is torsion free as well, so $\pi_*{\mathscr L}$ is torsion free. Together with 
$$ \dim_{\kappa(s)}\pi_*{\mathscr L} \otimes_{{\mathscr O}_S} \kappa(s) = 1 $$ 
for closed $s \in S$, this proves that $\pi_*{\mathscr L}$ is invertible. Indeed, the corresponding stalk of $\pi_*{\mathscr L}$ at $s$ is generated by a single element as an ${\mathscr O}_{s,S}$-module by Nakayama's lemma and, since it is torsion free, it has to be free of rank one. Finally, since  $(\varphi_S)_s: V \otimes_{\mathbb C} \kappa(s) \to \pi_*{\mathscr L} \otimes_{{\mathscr O}_S} \kappa(s)$ is surjective for closed points $s \in S$, we can define a morphism of schemes $\psi_S:S \to {\mathbb P}V$ such that $\text{Ker}(\varphi_S) = \psi_S^* {\mathscr U}_{{\mathbb P}V}$, where ${\mathscr U}_{{\mathbb P}V}$ is the tautological subbundle of $V \otimes {\mathscr O}_{{\mathbb P}V}$. \end{proof}

\begin{cor}\label{4.3} The sheaf $(\psi_S \circ \pi)^* {\mathscr U}_{{\mathbb P}V}$, regarded as an ${\mathscr O}_C$-submodule of $\mathrm{H}^0({\mathscr O}(1)) \otimes {\mathscr O}_C$, is contained in the kernel of $\mathrm{H}^0(X,{\mathscr O}(1)) \otimes {\mathscr O}_C \to f^*({\mathscr O}(1))$.\end{cor}

\begin{proof} By the definition of $\varphi_S$, $\text{Ker}(\varphi_S)$ lies inside the kernel of the map $\text{H}^0({\mathscr O}(1)) \otimes {\mathscr O}_S \to \pi_*f^*{\mathscr O}(1)$. First, $\pi^*$ is left exact because $\pi$ is flat, so $\pi^* \text{Ker}(\varphi_S) = (\psi_S \circ \pi)^* {\mathscr U}_{{\mathbb P}V}$ can be regarded as an ${\mathscr O}_C$-submodule of $\text{H}^0({\mathscr O}(1)) \otimes {\mathscr O}_C$. Moreover, it is contained in the kernel of the map 
$$ \text{H}^0({\mathscr O}(1)) \otimes {\mathscr O}_C \longrightarrow \pi^*\pi_*f^*({\mathscr O}(1)). $$
Composing the last map with $\pi^*\pi_*f^*{\mathscr O}(1) \to f^*{\mathscr O}(1)$, we recover the obvious pullback homomorphism $\text{H}^0({\mathscr O}(1)) \otimes {\mathscr O}_C \to f^*{\mathscr O}(1)$. \end{proof}

Next, we lift the map $\psi_S$ constructed above to the total space of the family of curves, $C$. The idea is obviously that $C_s$ is mapped by $f_s$ into $X_{\psi_S(s)}$, the vanishing locus of the codimension one subspace of $V \subset \text{H}^0({\mathscr O}(1))$ corresponding to $\psi_S(s) \in {\mathbb P}V$, but the formal statement will be given later in proposition \hyperlink{Prop4.5}{4.5}. We introduce $W$, the blowup of $X$ along $E$. As it is always the case with blowups, the graded ${\mathscr O}_X$-algebras homomorphism 
$$ \text{Sym} (V) \otimes {\mathscr O}_X \longrightarrow \bigoplus_{k=0}^\infty {\mathscr I}_{E/X}^k $$
induces the natural morphism of schemes $\tau = \tau_X \times \tau_{{\mathbb P} V}:W \to X \times {\mathbb P} V$. Note that there is a natural map $w:\text{proj}_{{\mathbb P}V}^* {\mathscr U}_{{\mathbb P}V} \to \text{proj}^*_X {\mathscr O}_X(1)$, where $\text{proj}_X$ and $\text{proj}_{{\mathbb P}V}$ denote the projections to the respective factors of $X \times {\mathbb P}V$.

\begin{claim}\label{4.4} The scheme-theoretic vanishing locus of $w$, regarded as a section of the vector bundle ${\mathscr H}om(\text{proj}_{{\mathbb P}V}^* {\mathscr U}_{{\mathbb P}V}, \text{proj}^*_X {\mathscr O}_X(1)) = {\mathscr O}_X(1) \boxtimes {\mathscr U}_{{\mathbb P}V}^\vee$, is precisely $W$.\end{claim}

\begin{proof} Let $\chi:X \to {\mathbb P}\text{H}^0(X,{\mathscr O}_X)$ be the embedding associated with ${\mathscr O}_X(1)$ and $P$ the blowup of ${\mathbb P}\text{H}^0(X,{\mathscr O}_X(1))$ along the projectivization of the cokernel of 
$$ V = \text{H}^0(X,{\mathscr I}_{E/X}(1)) \to \text{H}^0(X,{\mathscr O}_X(1)). $$
Again, $P$ sits naturally inside ${\mathbb P}\text{H}^0(X,{\mathscr O}_X(1)) \times {\mathbb P}V.$ Replacing $X$ with ${\mathbb P}\text{H}^0(X,{\mathscr O}_X)$, there is an analoguous way to define a natural section $p$ of ${\mathscr O}(1) \boxtimes {\mathscr U}_{{\mathbb P}V}^\vee$. The analoguous statement, that the vanishing locus of $p$ is $P$, is clear. Since $w$ is the pullback of $p$ to $X \times {\mathbb P}V$, the claim follows simply by pulling back via $\chi \times \text{Id}_{{\mathbb P}V}$. \end{proof}

\begin{prop}\label{4.5} There exists a lift $\tilde{f}:C \to W$ of $f$ along $\tau_X$ such that $\psi_S \circ \pi = \tau_{{\mathbb P} V} \circ \tilde{f}$, i.e. we require
\begin{center}
\begin{tikzpicture}[scale=0.8]

\matrix [column sep  = 15mm, row sep = 10mm] {
	& \node (nc) {$X$}; & \\
	\node (cw) {$C$}; & 
	\node (cc) {$W$}; & 
	\node (ce) {$X \times {\mathbb P}V$}; \\
	\node (sw) {$S$}; &
	\node (sc) {${\mathbb P}V$}; & \\
};

\draw[-{Stealth[length=2mm]}, thin] (cw) -- (nc);
\draw[-{Stealth[length=2mm]}, thin] (cw) -- (cc);
\draw[-{Stealth[length=2mm]}, thin] (cw) -- (sw);
\draw[-{Stealth[length=2mm]}, thin] (ce) -- (nc);
\draw[-{Stealth[length=2mm]}, thin] (cc) -- (ce);
\draw[-{Stealth[length=2mm]}, thin] (sw) -- (sc);
\draw[-{Stealth[length=2mm]}, thin] (ce) -- (sc);
\draw[-{Stealth[length=2mm]}, thin] (cc) -- (sc);
\draw[-{Stealth[length=2mm]}, thin] (cc) -- (nc);

\node at (-2,0.4) {$\tilde{f}$};
\node at (0.8,0.4) {$\tau$};
\node at (-2,-1.5) {$\psi_S$};

\end{tikzpicture}
\end{center}
to be a commutative diagram. \end{prop}

\begin{proof} If we show that the morphism $f \times (\psi_S \circ \pi): C \to X \times {\mathbb P} V$ factors through $\tau$, we may take $\tilde{f}$ to be the quotient morphism. The important point is that pullback by $f \times (\psi_S \circ \pi)$ kills $w$. Indeed, the pullback map $(f \times (\psi_S \circ \pi))^*w:(\psi_S \circ \pi)^* {\mathscr U}_{{\mathbb P}V} \to f^*{\mathscr O}_X(1)$ is the restriction to $(\psi_S \circ \pi)^* {\mathscr U}_{{\mathbb P}V}$ of the homomorphism $\text{H}^0(X,{\mathscr O}_X(1)) \otimes {\mathscr O}_C \to f^*{\mathscr O}_X(1)$, which is zero, by corollary \ref{4.3}. Then \ref{4.4} shows that $f \times (\psi_S \circ \pi): C \to X \times {\mathbb P} V$ indeed factors through $\tau$, completing the proof. \end{proof}

The crucial step is to understand the image of $\Phi_S$. The idea is very simple: if some $W_t$ contains some rational curve of some class $\smash{\tilde{\beta}}$ (lift of $\hat{\beta}$) through all $\xi_i$, then the rational curve will cut the copy of $E$ inside $W_t$ precisely in the divisor $\xi_1 + ... + \xi_e$ by degree considerations. The existence of divisor classes on $W_t$ restricting to a predetermined divisor class on $E$ imposes one condition on $t$.

Let $(C,\pi,f,p_1,p_2,...,p_e) \in \overline{\mathcal M}_e(X,E;\beta)|_{\Delta^c}(S)$ such that the map $\pi$ is smooth, and $W_S := W \times_S {\mathbb P}V \to S$ is smooth. The square in the diagram of Proposition \ref{4.5} induces an $S$-morphism $\smash{\hat{f}:C \to W_S}$. Note that, since we're over the complement of $\Delta$, any component of the source of any individual map is either contracted, or mapped birationally onto its image. Let ${\mathscr K}$ be the kernel of $\smash{ {\mathscr O}_{W_S} \to \hat{f}_*{\mathscr O}_C}$. The fact that $C$ and $W_S$ are flat over $S$ easily implies that ${\mathscr K}$ is also flat over $S$. However, since the restriction of ${\mathscr K}$ to $W_s$ for closed points $s \in S$ is invertible by obvious geometric considerations, it follows by \ref{2.4} that ${\mathscr K}$ is itself invertible. The section of ${\mathscr H}om({\mathscr K},{\mathscr O}_{W_S})$ corresponding to the inclusion ${\mathscr K} \to {\mathscr O}_{W_S}$ restricts on $E \times S$ to the section of ${\mathscr H}om({\mathscr I}_{p(S)},{\mathscr O}_{E \times S}) \cong {\mathscr O}_{E \times S}(p(S))$ corresponding to the inclusion ${\mathscr I}_{p(S)} \to {\mathscr O}_{E \times S}$, so we have a commutative square with $S$, $\text{Pic}(W_S/S)$, $E^e$ and $\text{Pic}(E)$.

By well-known properties, $\mathrm{H}_2(W,{\mathbf Z}) \cong \mathrm{H}_2(X,{\mathbf Z}) \oplus {\mathbf Z}$. Moreover, there exists a unique $(\tau_X)_*$-lift, $\hat{\beta} \in \text{H}_2(W,{\mathbf Z})$ of $\beta$, such that $(\tau_{{\mathbb P}V})_*\hat{\beta}= 0$. It is not hard to see that the mapping defined in Proposition \ref{4.5} $(C,\pi,f,p_1,p_2,...,p_e) \mapsto (C,\pi,\tilde{f},p_1,p_2,...,p_e)$ induces a morphism of Deligne-Mumford stacks
$$ \Phi:\overline{\mathcal M}_e(X,E;\beta)|_{\Delta^c} \to \overline{\mathcal M}_e(W,E \times {\mathbb P}V;\hat{\beta}). \eqno(4.3) $$
Since the obvious map $j:\overline{\mathcal M}_e(W,E \times {\mathbb P}V;\hat{\beta}) \to \overline{\mathcal M}_e(X,E;\beta)$ is inverse to $\Phi$, the morphism $\Phi$ induces an isomorphism $\overline{\mathcal M}_e(X,E;\beta)|_{\Delta^c} \to \overline{\mathcal M}_e(W,E \times {\mathbb P}V;\hat{\beta})|_{\Delta^c \times ({\mathbb P}V)^e}$. To conclude, there is a commutative diagram
\begin{center}
\begin{tikzpicture}

\matrix [column sep  = 15mm, row sep = 17mm] {
	\node (nw) {$\overline{\mathcal C}_e(X,E;\beta)|_{\Delta^c}$}; & 
	\node (nc) {$\overline{\mathcal C}_e(W,E \times {\mathbb P}V;\hat{\beta})$}; &
	\node (ne) {$W$}; \\
	\node (sw) {$\overline{\mathcal M}_e(X,E;\beta)|_{\Delta^c}$}; & 
	\node (sc) {$\overline{\mathcal M}_e(W,E \times {\mathbb P}V;\hat{\beta})$}; &
	\node (se) {${\mathbb P}V$}; \\
	\node (ssw) {$E^e$}; &
	\node (ssc) {$(E \times {\mathbb P}V)^e$}; & \\
};

\draw[-{Stealth[length=2mm]}, thin] (nw) -- (nc);
\draw[-{Stealth[length=2mm]}, thin] (nc) -- (ne);
\draw[-{Stealth[length=2mm]}, thin] (sw) -- (sc);
\draw[-{Stealth[length=2mm]}, thin] (sc) -- (se);
\draw[-{Stealth[length=2mm]}, thin] (nw) -- (sw);
\draw[-{Stealth[length=2mm]}, thin] (nc) -- (sc);
\draw[-{Stealth[length=2mm]}, thin] (ne) -- (se);
\draw[-{Stealth[length=2mm]}, thin] (sw) -- (ssw);
\draw[-{Stealth[length=2mm]}, thin] (sc) -- (ssc);
\draw[-{Stealth[length=2mm]}, thin] (ssc) -- (ssw);

\node at (-1.5,0.2) {$\Phi$};
\node at (3.3, 0.2) {$\Psi$};
\node at (-1.5, 2.6) {$\Phi_C$};
\node at (3.3, 2.6) {$\tilde{F}$};
\node at (-2.8, -1) {$\text{ev}_{X,E}$};
\node at (2, -1) {$\text{ev}_{W,E \times {\mathbb P}V}$};

\end{tikzpicture}
\end{center}
and $\Phi$ is an open embedding. Let $\xi \in E^e$ be a general $e$-tuple of closed points on $E$, so $\xi \in \Delta^c$. The diagram above shows that
$$ \text{ev}_X^{-1}(\xi) = \text{ev}_{X,E}^{-1}(\xi) \cong \text{ev}_{W,E \times {\mathbb P}V}^{-1}(\xi \times {\mathbb P}V). $$
Let ${\mathcal M}_\xi$ denote the last space and ${\mathcal C}_\xi$ its universal family, constructed in the obvious way by restricting the universal family above. 

Let $U \subset {\mathbb P}V$ be an open subset over which $W \to {\mathbb P}V$ is smooth. By a slight abuse of notation we write $\smash{ \overline{\mathcal M}_e(W,E \times {\mathbb P}V;\hat{\beta})|_U }$ where we actually mean restriction to $\Psi^{-1}(U) \cap (\Delta^c \times {\mathbb P}V^e)$. Denote by $\smash{ \overline{\mathcal M}^\circ_e(W,E \times {\mathbb P}V;\hat{\beta})|_U }$ the open locus where the source curve is smooth. By the discussion above, there is a map 
$$ \smash{ \omega: \overline{\mathcal M}^\circ_e(W,E \times {\mathbb P}V;\hat{\beta})|_U \to \text{Pic}(W_U/U)}, $$
where $W_U = \tau_{{\mathbb P}V}^{-1}(U)$. Let $\rho$ be the restriction map $\rho:\text{Pic}(W_U/U) \to U \times \text{Pic}(E)$. Similar to the construction in section 3, let $U'$ parametrize pairs $(t,\varphi_t)$, where $t \in U$ is a closed point and $\varphi_t$ is a ${\mathfrak D}$-isomorphism of $(\text{H}^2(W_t,{\mathbf Z}),\cup,c_1({\mathscr T}^\vee_{W_t}))$ with $\smash{(H_r,\langle \cdot , \cdot \rangle,\omega)}$. Then $U' \to U$ is an \'{e}tale morphism. The data is summarized in the following diagram.
\begin{center}
\begin{tikzpicture}

\matrix [column sep  = 15mm, row sep = 13mm] {
	\node (nw) {$U' \times_U \overline{\mathcal M}^\circ_e(W,E \times {\mathbb P}V;\hat{\beta})|_U$}; & 
	\node (nc) {$\text{Pic}(W_{U'}/U')$}; & \\
	\node (sw) {$\overline{\mathcal M}^\circ_e(W,E \times {\mathbb P}V;\hat{\beta})|_U$}; & 
	\node (sc) {$\text{Pic}(W_U/U)$}; &
	\node (se) {$\text{Pic}(E) \times U$}; \\ &
	\node (ssc) {$E^e$}; &
	\node (sse) {$\text{Pic}(E)$}; \\
};

\draw[-{Stealth[length=2mm]}, thin] (nw) -- (nc);
\draw[-{Stealth[length=2mm]}, thin] (sw) -- (sc);
\draw[-{Stealth[length=2mm]}, thin] (sc) -- (se);
\draw[-{Stealth[length=2mm]}, thin] (se) -- (sse);
\draw[-{Stealth[length=2mm]}, thin] (sw) -- (ssc);
\draw[-{Stealth[length=2mm]}, thin] (nw) -- (sw);
\draw[-{Stealth[length=2mm]}, thin] (nc) -- (sc);
\draw[-{Stealth[length=2mm]}, thin] (ssc) -- (sse);

\node at (-1,0.3) {$\omega$};
\node at (3.4, 0.3) {$\rho$};

\end{tikzpicture}
\end{center}
Finally, we remark that the analogue of Lemma \ref{3.1}, namely the statement that the section $U' \to \text{Pic}(W_{U'}/U')$ corresponding to an element $\beta \in H_r$ which is not a multiple of $\omega \in H_r$ composed with the map $\text{Pic}(W_{U'}/U') \to \text{Pic}(E)$ is nonconstant, follows from 3.1 itself simply by restricting to a line ${\mathbb P}^1 \subset {\mathbb P}V$. 

\begin{proof}[Proof of Theorem \ref{1.1}.] We will use the current notation rather than the notation $\smash{\mathrm{dP}_n^d}$ used in the statement of the theorems. Let $\smash{(C,{\tilde{f}},p_1,...,p_e) \in {\mathcal M}_\xi({\mathbf C})}$, which is the $\Phi$-image of some stable map $(C,f,p_1,...,p_e) \in \text{ev}_{X,E}^{-1}(\xi) \subset \overline{\mathcal M}_e(X,E;\beta) \subset \overline{\mathcal M}_e(X;\beta)$. Let
$$ {\mathscr N}_{f,X} = {\mathscr O}(a_1) \oplus {\mathscr O}(a_2) \oplus ... \oplus {\mathscr O}(a_{n-1}) $$
with $a_1 \geq a_2 \geq ... \geq a_{n-1}$. Then, as we've seen, $f$ maps to the surface $W_t \hookrightarrow X$. We denote the tangent space to ${\mathcal M}_\xi$ at $\smash{(C,{\tilde{f}},p_1,...,p_e)}$ by $\smash{\text{Def}^1_{W,\xi}({\tilde{f}})}$. Similarly, $\smash{\text{Def}^1_{W_t,\xi}({\tilde{f}})}$ denotes the space of first order deformations of $\smash{\tilde{f}}$ which remain inside $W_t$ and $\smash{\text{Def}^1_{X,\xi}(f)}$ denotes the space of first order deformations of $f$, subject to the condition that the marked points map to $\xi$. By the observation preceding the proof and the assumption that $W_t$ is general, we have an exact sequence
$$ 0 \longrightarrow \text{Def}^1_{W_t,\xi}({\tilde{f}}) \longrightarrow \text{Def}^1_{W,\xi}({\tilde{f}}) \longrightarrow {\mathscr T}_t {\mathbb P}V \longrightarrow {\mathbf C} \longrightarrow 0. $$
It follows that $\smash{ \dim \text{Def}^1_{W,\xi}({\tilde{f}}) = \dim \text{Def}^1_{W_t,\xi}({\tilde{f}}) + \dim {\mathbb P}V - 1}$. We claim that
$$\text{Def}^1_{W_t,\xi}({\tilde{f}}) = 0. $$
Recall that $\smash{\tilde{f}}$ is unramified, so the normal bundle $\smash{{\mathscr N}_{\tilde{f},W_t}}$ of $\smash{\tilde{f}}$ relative to $W_t$ is locally free of rank one. From the standard sequence
$$ 0 \longrightarrow {\mathscr T}_C \longrightarrow \tilde{f}^*{\mathscr T}_{W_t} \longrightarrow {\mathscr N}_{\tilde{f},W_t} \longrightarrow 0, $$
we infer that $\smash{ c_1({\mathscr N}_{\tilde{f},W_t}) =\tilde{f}^*c_1(W_t) - c_1({\mathscr T}_C)}$, so $\smash {\deg {\mathscr N}_{\tilde{f},W_t} =e-2}$. Therefore,
$$ \text{Def}^1_{W_t,\xi}({\tilde{f}}) = \text{H}^0({\mathscr N}_{\tilde{f},W_t} \otimes {\mathscr O}_C(-p_1-...-p_e)) = 0, $$
as desired. Therefore, $\smash{ \dim \text{Def}^1_{W,\xi}({\tilde{f}}) = n-3}$, so $\smash{\text{Def}^1_{X,\xi}(f) = n-3}$ as well, since $\Phi$ is an open immersion. Using once more the interplay between first order deformations and normal bundles, we conclude that $h^0({\mathscr N}_{f,X}(-p_1-...-p_e)) = n-3$. Therefore,
$$ \dim \bigoplus_{i=1}^{n-1} \text{H}^0({\mathscr O}(a_i-e)) = n-3, $$
and combining with the constrain $a_1+a_2+...+a_{n-1} = \deg c_1({\mathscr N}_{f,X}) = (n-1)e-2$, we can infer that $a_{n-1} \geq e-1$.

Let $\xi_{e+1} \in X \backslash E = W \backslash E \times {\mathbb P}V$, say $\xi_{e+1} \in W_t$. We consider stable maps with one additional marked point. Let $\xi_1,...,\xi_e$ as before, $\xi'=(\xi_1,...,\xi_e,\xi_{e+1})$ and $\text{ev}_X^{-1}(\xi')$ the space of stable maps $f$ to $X$ such that $f(p_i) = \xi_i$. We will only sketch this part of the argument, being similar to the analysis above. The condition $f(p_{e+1}) = \xi_{e+1}$ forces $f$ to map to $W_t$, so $\smash{\text{Def}^1_{X,\xi'}(f)}$ is now isomorphic to $\smash{\text{Def}^1_{W_t,\xi'}({\tilde{f}})}$, but, as above, $\smash{\text{Def}^1_{W_t,\xi'}({\tilde{f}})} = 0$, so $\smash{\text{Def}^1_{X,\xi'}(f)} = 0$. In terms of normal bundles, $h^0({\mathscr N}_{f,X}(-p_1-...-p_{e+1})) = 0$, implying that $a_1 \leq e$. Combining with the previous inequality, we conclude that $a_1=a_2=e$ and $a_3=...=a_{n-1} = e-1$, as desired.  \end{proof}

\begin{proof}[Proof of Theorem \ref{1.2}.] We are in the situation analyzed in section 3, $X = \text{dP}_3^d$ a Fano threefold of index 2 and degree $d \in \{2,3,4,5\}$ and $W \to {\mathbb P}^1$ a Lefschetz pencil. Let ${\mathscr L}_E := {\mathscr O}_E(\xi_1+...+\xi_e) \in {\text{Pic}}^e(E)$. As in section 3, let $S(\xi) = S({\mathscr L}_E)$ be the set of pairs $(t,{\mathscr L}_t)$ such that $t \in {\mathbb P}^\circ = U$ and ${\mathscr L}_t$ restricts to ${\mathscr L}_E$ on $E$. To each such ${\mathscr L}_t$, we can associate $\smash{\tilde{\beta}_t} \in \text{H}_2(W_t,{\mathbf Z})$, the Poincar\'{e} dual to $c_1({\mathscr L}_t)$. The condition in the second part of Lemma \ref{2.1} is satisfied for a general choice of $\xi$. Indeed, for $e \geq 2$, we may move the $e$ points around preserving ${\mathscr L}_E$ (and therefore $S(\xi)$) to avoid the finitely many prohibited situations. For $e=1$, the condition is satisfied vacuously. 

Of course, $S(\xi)$ splits as a disjoint union of $S({\mathscr L}_E,\gamma)$ consisting of those pairs $(t,{\mathscr L}_t)$ such that $\smash{\tilde{\beta}_t}$ corresponds to $\gamma$ under a suitable ${\mathfrak D}$-isomorphism, where $\gamma$ varies over a set of representatives of $H^e_r/G$. From the discussion at the end of section 3, $S({\mathscr L}_E,\gamma)$ is reduced and we have
\begin{equation*}
\begin{aligned} {\mathcal M}_\xi &\cong \bigsqcup_{(t,{\mathscr L}_t) \in S(\xi)} \overline{\mathcal M}_e(W_t,\tilde{\beta}_t) \times_{W_t^e} \{ \xi \}\\
&\cong \bigsqcup_{(t,{\mathscr L}_t) \in S(\xi)} \overline{\mathcal M}_{e-1}(W_t,\tilde{\beta}_t) \times_{W_t^{e-1}} \{ (\xi_1,...,\xi_{e-1}) \}.
\end{aligned}
\end{equation*}
From the final observation in section 2.1, taking the degrees of these $0$-cycles, we obtain
$$ \left\langle [\text{pt}]^e \right\rangle_{0,\beta}^X = \sum_{(t,{\mathscr L}_t) \in S(\xi)} \left\langle [\text{pt}]^{e-1} \right\rangle_{0,\tilde{\beta}_t}^{W_t} = \sum_{\gamma \in H_r^e/G} \sum_{(t,{\mathscr L}_t) \in S({\mathscr L}_E,\gamma)} \left\langle [\text{pt}]^{e-1} \right\rangle_{0,\tilde{\beta}_t}^{W_t} $$
However, all $W_t$ are deformation equivalent to any degree $d$ del Pezzo surface $\Sigma$, so $\smash{  \left\langle [\text{pt}]^{e-1} \right\rangle_{0,\tilde{\beta}_t}^{W_t} }$ is $\smash{ \left\langle [\text{pt}]^{e-1} \right\rangle_{0,\gamma}^{\Sigma} }$ for all $(t,{\mathscr L}_t) \in S({\mathscr L}_E,\gamma)$. By (3.5), we conclude that
$$ \left\langle [\text{pt}]^e \right\rangle_{0,\beta}^X = \sum_{\gamma \in H_r^e/G} |S({\mathscr L}_E,\gamma)| \left\langle [\text{pt}]^{e-1} \right\rangle_{0,\gamma}^{\Sigma} = \frac{1}{d(9-d)} \sum_{\gamma \in H^e_r}\Delta(\omega,\gamma)\left\langle [\text{pt}]^{e-1} \right\rangle_{0,\gamma}^{\Sigma}, $$
completing the proof of the formula.\end{proof}

\end{document}